\documentclass[11pt]{article}
\usepackage{amsmath}
\usepackage{amssymb}
\usepackage{latexsym}
\begin{document}
\title{A twisted generalization of Lie-Yamaguti algebras}
\author{Donatien Gaparayi \\ Institut de Math\'ematiques et de Sciences Physiques, \\
01 BP 613, Porto-Novo, BENIN \\ and \\ A. Nourou Issa \\ D\'epartement de Math\'ematiques, Universit\'e d'Abomey-Calavi, \\ 01 BP 4521, Cotonou 01, BENIN. {\it Email}: woraniss@yahoo.fr}
\date{}
\maketitle
\begin{abstract}
A twisted generalization of Lie-Yamaguti algebras, called Hom-Lie-Yamaguti algebras, is defined. Hom-Lie-Yamaguti algebras generalize Hom-Lie triple systems (and subsequently ternary Hom-Nambu algebras) and Hom-Lie algebras in the same way as Lie-Yamaguti algebras generalize Lie triple systems and Lie algebras. It is shown that the category of Hom-Lie-Yamaguti algebras is closed under twisting by self-morphisms. Constructions of Hom-Lie-Yamaguti algebras from classical Lie-Yamaguti algebras and Malcev algebras are given. It is observed that, when the ternary operation of a Hom-Lie-Yamaguti algebra expresses through its binary one in a specific way, then such a Hom-Lie-Yamaguti algebra is a Hom-Malcev algebra.
\end{abstract}
\section{Introduction}
Using the Bianchi identities, K. Nomizu [14] characterized, by some identities involving the torsion and the curvature, reductive homogeneous spaces with canonical connection. K. Yamaguti [17] gave an algebraic interpretation of these identities by considering the torsion and curvature tensors of Nomizu's canonical connection as a bilinear and a trilinear algebraic operations satisfying some axioms, and thus defined what he called a ``general Lie triple system''. M. Kikkawa [7] used the term ``Lie triple algebra'' for such an algebraic object. More recently, M.K. Kinyon and A. Weinstein [8] introduced the term ``Lie-Yamaguti algebra'' for this object.
\par
A {\it Lie-Yamaguti algebra} $(V,*, \{ , , \})$ is a vector space $V$ together with a binary operation $*: V \times V \rightarrow V$ and a ternary operation $\{ , , \} : V \times V \times V \rightarrow V$ such that \\
\par
{\bf (A1)} $x*y = - y*x$,
\par
{\bf (A2)} $\{x,y,z\} = - \{y,x,z\}$,
\par
{\bf (A3)} $\sigma [(x*y) * z + \{x,y,z\}] = 0$,
\par
{\bf (A4)} $\sigma [\{ x*y , z, u \} = 0$,
\par
{\bf (A5)} $\{x, y, u*v \} = \{x,y,u\}*v + u*\{x,y,v\}$,
\par
{\bf (A6)} $ \{x,y, \{u,v,w \} \} = \{ \{ x,y,u \},v,w \} + \{u, \{ x,y,v \},w \}$
\par
\hspace{4cm} $+ \{ u,v, \{x,y,w\} \}$, \\
for all $u,v,w,x,y,z$ in $V$ and $\sigma$ denotes the sum over cyclic permutation of $x,y,z$. \\
\par
In [2] the notation ``LY-algebra'' is used for ``Lie-Yamaguti algebra''. So, likewise, we will write ``Hom-LY algebra'' for ``Hom-Lie-Yamaguti algebra''.
\par
Observe that if $x*y = 0$, for all $x,y$ in $V$, then $(V,*, \{ , , \})$ reduces to a {\it Lie triple system} $(V,\{ , , \})$ as defined in [16]. Originally, N. Jacobson [6] defined a Lie triple system as a submodule of an associative algebra that is closed under the iterated commutator bracket.
\par
In this paper we consider a Hom-type generalization of LY algebras that we call Hom-LY algebras. Roughly, a Hom-type generalization of a given type of algebras is defined by twisting the defining identities of that type of algebras by a self-map in such a way that, when the twisting map is the identity map, one recovers the original type of algebras. The systematic study of Hom-algebras was initiated by A. Makhlouf and S.D. Silvestrov [11], while D. Yau [20] gave a general construction method of Hom-type algebras starting from classical algebras and a twisting self-map. For information on various types of Hom-algebras, one may refer to [1], [4], [5], [9]-[11], [19]-[22].
\par
A Hom-type generalization of $n$-ary Lie algebras, $n$-ary Nambu algebras and $n$-ary Nambu-Lie algebras (i.e. Fillipov $n$-ary algebras) called $n$-ary Hom-Lie algebras, $n$-ary Hom-Nambu algebras and $n$-ary Hom-Nambu-Lie algebras respectively, is considered in [1]. Such a generalization  is extended to the one of Hom-Lie triple systems and Hom-Jordan triple systems in [22]. We point out that the class of Hom-LY algebras encompasses the ones of ternary Hom-Nambu algebras, Hom-Lie triple systems (hence Jordan and Lie triple systems), Hom-Lie algebras (hence Lie algebras) and LY algebras.
\par
The rest of the paper is organized as follows. In section 2 some basic facts on Hom-algebras and $n$-ary Hom-algebras are recalled. The emphasis point here is that the definition of a Hom-triple system (Definition 2.3) is more restrictive than the D. Yau's in [22]. However, with this vision of a Hom-triple system, we point out that any non-Hom-associative algebra (i.e. nonassociative Hom-algebra or Hom-nonassociative algebra) has a natural structure of a Hom-triple system (this is the Hom-counterpart of a similar well-known result connecting nonassociative algebras and triple systems). Then we give the definition of a Hom-LY algebra and make some observations on its relationships with some types of ternary Hom-algebras and with LY algebras. In section 3 we show that the category of Hom-LY algebras is closed under twisting by self-morphisms (Theorem 3.1). Subsequently, we show a way to construct Hom-LY algebras from LY algebras (or Malcev algebras) by twisting along self-morphisms (Corollary 3.2 and Corollary 3.3); this is an extension to binary-ternary algebras of a result due to D. Yau ([20], Theorem 2.3. Such an extension is first mentioned in [4], Corollary 4.6). In section 4 some relationships between Hom-LY algebras and Hom-Malcev algebras are considered. We show that when the ternary operation of a Hom-LY algebra expresses through its binary one in a specific way, then such a Hom-LY algebra turns out to be a Hom-Malcev algebra (Proposition 4.1). Moreover, with this expression of the ternary operation, it is observed that, in a Hom-Malcev algebra, the Hom-Malcev identity can be written in terms of this ternary operation and the original binary operation of the given Hom-Malcev algebra (Proposition 4.2). These considerations constitute the Hom-version of similar relationships between Malcev algebras and LY algebras ([13], [18]).
\par
All vector spaces and algebras throughout will be over a ground field $\mathbb K$ of characteristic $0$.

\section{Ternary Hom-algebras. Definitions}
We recall some basic facts about Hom-algebras, including ternary Hom-Nambu algebras. We note that the definition of a Hom-triple system given here (see Definition 2.3) is slightly more restrictive than the one given by D. Yau [22]. Then we give the definition of the main object of this paper (see Definition 2.6) and show its relationships with known structures such as ternary Hom-Nambu algebras, Hom-Lie triple systems, Hom-Lie algebras or Lie-Yamaguti algebras.
\par
For definitions of $n$-ary Hom-algebras ($n$-ary Hom-Nambu and Hom-Nambu-Lie algebras, $n$-ary Hom-Lie algebras, etc.) we refer to [1], [22]. Here, for our purpose, we restrict our concern to ternary Hom-algebras. In fact, as we shall see below, a Hom-Lie-Yamaguti algebra is a ternary Hom-Nambu algebra with an additional binary anticommutative operation satisfying some compatibility conditions. \\
\par
{\bf Definition 2.1.} [22] A {\it ternary Hom-algebra} $(V, [ , , ], \alpha = ({\alpha}_{1}, {\alpha}_{2}))$ consists of a $\mathbb K$-module $V$, a trilinear map $[ , , ] : V \times V \times V \rightarrow V$, and linear maps ${\alpha}_{i} : V \rightarrow V$, $i = 1,2$, called the {\it twisting maps}. The algebra $(V, [ , , ], \alpha = ({\alpha}_{1}, {\alpha}_{2}))$ is said {\it multiplicative} if ${\alpha}_{1} = {\alpha}_{2} = \alpha $ and $\alpha ([x,y,z]) = [\alpha (x), \alpha (y), \alpha (z)]$ for all $x,y,z \in V$. \\
\par
For convenience, we assume throughout this paper that all Hom-algebras are multiplicative. \\
\par
{\bf Definition 2.2.} [1] A {\it (multiplicative) ternary Hom-Nambu algebra} is a (multiplicative) ternary Hom-algebra $(V, [ , , ], \alpha)$ satisfying \\
\par
$[\alpha (x), \alpha (y), [u,v,w]] = [[x,y,u], \alpha (v), \alpha (w)] + [\alpha (u), [x,y,v], \alpha (w)]$
\par
\hspace{3.5cm} $+ [\alpha (u), \alpha (v), [x,y,w]]$, \hfill (2.1) \\
\\
for all $u,v,w,x,y \in V$. \\
\par
The condition (2.1) is called the {\it ternary Hom-Nambu identity}. \\
\par
{\bf Definition 2.3.} A {\it (multiplicative) Hom-triple system} is a (multiplicative) ternary Hom-algebra $(V, [ , , ], \alpha)$ such that
\par
(i) $[u,v,w] = -[v,u,w]$,
\par
(i) $\sigma [u,v,w] = 0$, \\
for all $u,v,w \in V$, where $\sigma [u,v,w] := [u,v,w] + [v,w,u] + [w,u,v]$. \\
\par
{\bf Remark.} A more general definition of a Hom-triple system is given by D. Yau [22] without the requirements (i), (ii) as in Definition 2.3 above. Our definition here is motivated by the concern of giving a Hom-type analogue of the relationships between nonassociative algebras and triple systems (see Remark below). \\
\par
A Hom-algebra in which the Hom-associativity is not assumed is called a nonassociative Hom-algebra [10] or a Hom-nonassociative algebra [19] (the expression of ``non-Hom-associative'' Hom-algebra is used in [4] for that type of Hom-algebras). With the notion of a Hom-triple system as above, we have the following \\
\par
{\bf Proposition 2.4.} {\it Any non-Hom-associative Hom-algebra is a Hom-triple system}. \\
\par
{\bf Proof.} Let $(A, \cdot , \alpha)$ be a non-Hom-associative algebra. Then \\ 
$(A, [,], as( , , ), \alpha)$ is a Hom-Akivis algebra with respect to $[x,y] := x \cdot y - y \cdot x$ (commutator) and $as(x,y,z) := xy \cdot \alpha (z) - \alpha (x) \cdot yz$ (Hom-associator), i.e. the Hom-Akivis identity
\par
$\sigma [[x,y], \alpha (z)] = \sigma as(x,y,z) - \sigma as(y,x,z)$ \\
holds for all $x,y,z$ in $A$ ([4]). Now define
\par
$[x,y,z] := [[x,y], \alpha (z)] - as(x,y,z) + as(y,x,z)$ \\
for all $x,y,z$ in $A$. Then $[x,y,z] = -[y,x,z]$ and the Hom-Akivis identity implies that $\sigma [x,y,z] = 0$. Thus $(A, [, ,] , \alpha)$ is a Hom-triple system. \hfill $\square$ \\
\par
{\bf Remark.} For $\alpha = Id$ (the identity map), we recover the triple system with ternary operation $[[x,y],z] - (x,y,z) + (y,x,z)$ that is associated to each nonassociative algebra, since any nonassociative algebra has a natural Akivis algebra structure with respect to the commutator and associator operations $[x,y]$ and $(x,y,z)$, for all $x,y,z$ (see, e.g., remarks in [4]). \\
\par
{\bf Definition 2.5.} [22] A {\it Hom-Lie triple system} is a Hom-triple system $(V, [, ,] , \alpha)$ satisfying the ternary Hom-Nambu identity (2.1). \\
\par
When $\alpha = Id$, a Hom-Lie triple system reduces to a Lie triple system.
\par
We now give the definition of the basic object of this paper. \\
\par
{\bf Definition 2.6.} A {\it Hom-Lie-Yamaguti algebra} (Hom-LY algebra for short) is a quadruple $(L,*, \{ , , \}, \alpha)$ in which $L$ is a $\mathbb K$-vector space, ``$*$'' a binary operation and ``$\{ , , \}$'' a ternary operation on $L$, and $\alpha : L \rightarrow L$ a Linear map such that
\par
{\bf (B1)} $\alpha ( x * y) = \alpha (x) * \alpha (y)$,
\par
{\bf (B2)} $\alpha (\{x,y,z\}) = \{\alpha (x), \alpha (y), \alpha (z) \}$,
\par
{\bf (B3)} $x*y = - y*x$,
\par
{\bf (B4)} $\{x,y,z\} = - \{y,x,z\}$,
\par
{\bf (B5)} $\sigma [(x*y) * \alpha (z) + \{x,y,z\}] = 0$,
\par
{\bf (B6)} $\sigma [\{ x*y , \alpha (z), \alpha (u) \} + \{ z*y , \alpha (x), \alpha (u) \}] = 0$,
\par
{\bf (B7)} $\{\alpha (x), \alpha (y), u*v \} = \{x,y,u\}* {\alpha}^{2}(v) + {\alpha}^{2}(u)*\{x,y,v\}$,
\par
{\bf (B8)} $ \{ {\alpha}^{2}(x), {\alpha}^{2}(y), \{ {\alpha}^{2}(u), {\alpha}^{2}(v), w \} \} = \{ \{ x,y, {\alpha}^{2}(u) \}, {\alpha}^{4}(v), {\alpha}^{2}(w) \}$
\par
\hspace{6.5cm} $+ \{ {\alpha}^{4}(u), \{ x,y, {\alpha}^{2}(v) \}, {\alpha}^{2}(w) \}$
\par
\hspace{6.5cm} $+ \{ {\alpha}^{4}(u), {\alpha}^{4}(v), \{x,y,w\} \}$, \\
for all $u,v,w,x,y,z$ in $L$ and $\sigma$ denotes the sum over cyclic permutation of $x,y,z$. \\
\par
Note that the conditions {\bf (B1)} and {\bf (B2)} mean the multiplicativity of $(L,*, \{ , , \}, \alpha)$. \\
\par
{\bf Remark.} (1) If $\alpha = Id$, then the Hom-LY algebra $(L,*, \{ , , \}, \alpha)$ reduces to a LY algebra $(L,*, \{ , , \})$ (see {\bf (A1)}-{\bf (A6)}).
\par
(2) If $x*y = 0$, for all $x,y \in L$, then $(L,*, \{ , , \}, \alpha)$ becomes a Hom-Lie triple system $(L,\{ , , \}, {\alpha}^{2})$ and, subsequently, a ternary Hom-Nambu algebra (since, by Definition 2.5, any Hom-Lie triple system is automatically a ternary Hom-Nambu algebra).
\par
(3) If $\{x,y,z\} = 0$ for all $x,y,z \in L$, then the Hom-LY algebra $(L,*, \{ , , \}, \alpha)$ becomes a Hom-Lie algebra $(L,*, \alpha)$.
\section{Constructions of Hom-Lie-Yamaguti algebras}
In this section we consider construction methods for Hom-LY algebras. These methods allow to find examples of Hom-LY algebras starting from classical LY algebras or even from Malcev algebras.
\par
First, as the main tool, we show that the category of (multiplicative) Hom-LY algebras is closed under self-morphisms. \\
\par
{\bf Theorem 3.1.} {\it Let $A_{\alpha} := (A,*,\{ , , \}, \alpha)$ be a multiplicative Hom-LY algebra and let $\beta$ be an endomorphism of the algebra $(A,*,\{ , , \})$ such that $\beta \alpha = \alpha \beta$. Define on $A$ the operations}
\par
$x *_{\beta} y := \beta (x * y)$,
\par
$\{x,y,z \}_{\beta} := {\beta}^{2}(\{x,y,z \})$ \\
{\it for all $x,y,z$ in $A$. Then $A_{\beta} := (A, *_{\beta}, \{ ,, \}_{\beta}, \beta \alpha)$ is a multiplicative Hom-LY algebra}. \\
\par
{\bf Proof.} We have
\par
$ (\beta \alpha)(x *_{\beta} y) = (\beta \alpha)(\beta (x) * \beta (y)) = \beta ((\alpha \beta)(x) * (\alpha \beta)(y)) $
\par
$= (\alpha \beta)(x) *_{\beta} (\alpha \beta)(y) = (\beta \alpha)(x) *_{\beta} (\beta \alpha)(y)$ and we get {\bf (B1)}. Likewise, the condition $\beta \alpha = \alpha \beta$ implies {\bf (B2)}. The identities {\bf (B3)} and {\bf (B4)} for $A_{\beta}$ follow from the skew-symmetry of ``$*$'' and ``$\{ , , \}$`` respectively.
\par
Consider now $\sigma ((x *_{\beta} y)  *_{\beta} (\beta \alpha)(z)) + \sigma \{x,y,z \}_{\beta}$. Then
\begin{eqnarray}
\sigma ((x *_{\beta} y)  &*_{\beta}& (\beta \alpha)(z)) + \sigma \{x,y,z \}_{\beta} \nonumber \\
&=& \sigma [\beta (\beta (x * y) * \beta (\alpha (z)))] + \sigma [{\beta}^{2} (\{x,y,z \})] \nonumber \\
&=& \sigma [\beta ((\beta (x) * \beta (y)) * \alpha (\beta (z))] + \sigma [{\beta}^{2} (\{x,y,z \})] \nonumber \\ 
& & \mbox{(since} \; \beta \alpha = \alpha \beta) \nonumber \\
&=& \beta (\sigma [(\beta (x) * \beta (y)) * \alpha (\beta (z))] + \sigma \{\beta (x), \beta (y), \beta (z) \}) \nonumber \\
&=& \beta (0) \; \mbox{(by {\bf (B6)} for} \; A_{\alpha}) \nonumber \\
&=& 0 \nonumber
\end{eqnarray}
and thus we get {\bf (B5)} for $A_{\beta}$. Next,
\par
$ \{ x *_{\beta} y, (\beta \alpha)(z), (\beta \alpha)(u) \}_{\beta} = \{ {\beta}^{3}(x*y), {\beta}^{3}(\alpha (z)), {\beta}^{3}(\alpha (u)) \}$
\par
$ = {\beta}^{3}(\{x*y, \alpha (z), \alpha (u) \})$. \\
Likewise we find that $ \{ z *_{\beta} y, (\beta \alpha)(x), (\beta \alpha)(u) \}_{\beta} = {\beta}^{3}(\{z*y, \alpha (x), \alpha (u) \})$. Therefore
\begin{eqnarray}
& & \{ x*_{\beta}y , (\beta \alpha)(z), (\beta \alpha)(u) \}_{\beta} + \{ z *_{\beta} y, (\beta \alpha)(x), (\beta \alpha)(u) \}_{\beta} \nonumber \\
&=& \sigma [{\beta}^{3}(\{x*y, \alpha (z), \alpha (u) \}) + {\beta}^{3}(\{z*y, \alpha (x), \alpha (u) \})] \nonumber \\
&=& {\beta}^{3}(\sigma \{x*y, \alpha (z), \alpha (u) \} + \sigma \{z*y, \alpha (x), \alpha (u) \}) \nonumber \\
&=& {\beta}^{3} (0) \; \mbox{(by {\bf (B6)} for} \; A_{\alpha}) \nonumber \\
&=& 0 \nonumber
\end{eqnarray}
so that we get {\bf (B6)} for $A_{\beta}$. Further, using {\bf (B7)} for $A_{\alpha}$ and condition $\alpha \beta = \beta \alpha$, we compute
\par
$\{(\beta \alpha)(x), (\beta \alpha)(y), u \; *_{\beta} \; v \}_{\beta} = {\beta}^{3}(\{ \alpha (x), \alpha (y), u*v \})$ \\ $= {\beta}^{3}( \{ x, y, u \} * {\alpha}^{2}(v) + {\alpha}^{2}(u) * \{ x, y, v \}) = \beta ({\beta}^{2}(\{ x, y, u \}) * ({\beta}^{2} {\alpha}^{2})(v))$ \\ $ + \beta (({\beta}^{2} {\alpha}^{2})(u) * {\beta}^{2}(\{ x, y, v \})) = \{ x, y, u \}_{\beta} \; *_{\beta} \; ({\beta}^{2} {\alpha}^{2})(v)$ \\ $ + ({\beta}^{2} {\alpha}^{2})(u) \; *_{\beta} \; \{ x, y, v \}_{\beta} $ \\ $= \{ x, y, u \}_{\beta} \; *_{\beta} \; {(\beta \alpha)}^{2}(v) + {(\beta \alpha)}^{2}(u) \; *_{\beta} \; \{ x, y, v \}_{\beta}$. \\
Thus {\bf (B7)} holds for $A_{\beta}$. Using repeatedly the condition $\alpha \beta = \beta \alpha$ and the identity {\bf (B8)} for $A_{\alpha}$, the verification of {\bf (B8)} for $A_{\beta}$ is as follows.
\begin{eqnarray}
& & \{ {(\beta \alpha)}^{2}(x), {(\beta \alpha)}^{2}(y), \{ {(\beta \alpha)}^{2}(u), {(\beta \alpha)}^{2}(v), w \}_{\beta} \}_{\beta} \nonumber \\
&=& \{({\beta}^{2} {\alpha}^{2})(x), ({\beta}^{2} {\alpha}^{2})(y), \{({\beta}^{2} {\alpha}^{2})(u), ({\beta}^{2} {\alpha}^{2})(v), w\}_{\beta} \}_{\beta} \nonumber \\
&=&  {\beta}^{2}(\{ ({\beta}^{2} {\alpha}^{2})(x), ({\beta}^{2} {\alpha}^{2})(y), {\beta}^{2}(\{({\beta}^{2} {\alpha}^{2})(u), ({\beta}^{2} {\alpha}^{2})(v), w \})\}) \nonumber \\
&=& {\beta}^{4}(\{ {\alpha}^{2}(x), {\alpha}^{2}(y), \{({\beta}^{2} {\alpha}^{2})(u), ({\beta}^{2} {\alpha}^{2})(v), w \} \}) \nonumber \\
&=& {\beta}^{4}(\{ {\alpha}^{2}(x), {\alpha}^{2}(y), \{{\alpha}^{2}({\beta}^{2}(u)), {\alpha}^{2}({\beta}^{2}(v)), w \} \}) \nonumber \\
&=& {\beta}^{4}(\{{\alpha}^{4}({\beta}^{2}(u)), {\alpha}^{4}({\beta}^{2}(v)), \{x,y,w \} \}) \nonumber \\
&+& {\beta}^{4}(\{ \{x,y, {\alpha}^{2}({\beta}^{2}(u))\}, {\alpha}^{4}({\beta}^{2}(v)), {\alpha}^{2}(w) \} \nonumber \\
&+& {\beta}^{4}(\{{\alpha}^{4}({\beta}^{2}(u)), \{x,y, {\alpha}^{2}({\beta}^{2}(v))\}, {\alpha}^{2}(w) \} \nonumber \\
&=& {\beta}^{4}(\{({\beta}^{2} {\alpha}^{4})(u), ({\beta}^{2} {\alpha}^{4})(v), \{x,y,w \} \}) \nonumber \\
&+& {\beta}^{4}(\{ \{x,y,({\beta}^{2} {\alpha}^{2})(u)\}, ({\beta}^{2} {\alpha}^{4})(v), {\alpha}^{2}(w) \}) \nonumber \\
&+& {\beta}^{4}(\{({\beta}^{2} {\alpha}^{4})(u), \{x,y,({\beta}^{2} {\alpha}^{2})(v)\}, {\alpha}^{2}(w) \}) \nonumber \\
&=&  {\beta}^{2}(\{({\beta}^{4} {\alpha}^{4})(u), ({\beta}^{4} {\alpha}^{4})(v), {\beta}^{2}(\{x,y,w \})\}) \nonumber \\
&+& {\beta}^{2}(\{ {\beta}^{2}(\{x,y,({\beta}^{2} {\alpha}^{2})(u)\}), ({\beta}^{4} {\alpha}^{4})(v), ({\beta}^{2} {\alpha}^{2})(w)\}) \nonumber \\
&+& {\beta}^{2}(\{({\beta}^{4} {\alpha}^{4})(u),{\beta}^{2}(\{x,y,({\beta}^{2} {\alpha}^{2})(v)\}), ({\beta}^{2} {\alpha}^{2})(w)\}) \nonumber \\
&=& \{ {(\beta \alpha)}^{4}(u), {(\beta \alpha)}^{4}(v), \{x,y,w \}_{\beta} \}_{\beta} \nonumber \\
&+& \{ \{x,y, {(\beta \alpha)}^{2}(u) \}_{\beta}, {(\beta \alpha)}^{4}(v), {(\beta \alpha)}^{2}(w) \}_{\beta} \nonumber \\
&+& \{ {(\beta \alpha)}^{4}(u),\{x,y, {(\beta \alpha)}^{2}(v) \}_{\beta}, {(\beta \alpha)}^{2}(w) \}_{\beta}. \nonumber
\end{eqnarray}
Thus {\bf (B8)} holds for $A_{\beta}$. Therefore, we get that $A_{\beta}$ is a Hom-LY algebra. This finishes the proof. \hfill $\square$ \\
\par
From Theorem 3.1 we have the following method of construction of Hom-LY algebras from LY algebras. This method is an extension to binary-ternary algebras of a result due to D. Yau ([20], Theorem 2.3), giving a general method of construction of Hom-algebras from their corresponding untwisted algebras. Such an extension to binary-ternary algebras is first mentioned in [4], Corollary 4.6. \\
\par
{\bf Corollary 3.2.} {\it Let $(A, * , [ , , ])$ be a LY algebra and $\beta$ an endomorphism of $(A, * , [ , , ])$. If define on $A$ a binary operation ''$\tilde{*}$`` and a ternary operation ''$\{ , , \}$`` by}
\par
$x \tilde{*} y := \beta (x * y)$,
\par
$\{ x,y,z \} := {\beta}^{2}([x,y,z])$, \\
{\it then $(A, \tilde{*}, \{ , , \}, \beta)$ is a Hom-LY algebra}. \\
\par
{\bf Proof.} The proof follows if observe that Corollary 3.2 is Theorem 3.1 when $\alpha = Id$. \hfill $\square$ \\
\par
{\bf Corollary 3.3.} {\it Let $(A,*)$ be a Malcev algebra and $\beta$ any endomorphism of $(A,*)$. Define on $A$ the operations}
\par
$x \tilde{*} y := \beta (x * y)$,
\par
$\{ x,y,z \} := {\beta}^{2}((x*y)*z - (y*z)*x - (z*x)*y)$. \\
{\it Then $(A, \tilde{*}, \{ , , \}, \beta)$ is a Hom-LY algebra}. \\
\par
{\bf Proof.} If consider on $A$ the ternary operation $[x,y,z] :=$ \\ $(x*y)*z - (y*z)*x - (z*x)*y$, $\forall x,y,z \in A$, then $(A, * , [ , , ])$ is a LY algebra [18]. Moreover, since $\beta$ is an endomorphism of $(A,*)$, we have $\beta ([x,y,z]) = (\beta(x) * \beta(y)) * \beta(z) - (\beta(y) * \beta(z)) * \beta(x) - (\beta(z) * \beta(x)) * \beta(y) = [\beta(x),\beta(y),\beta(z)]$ so that $\beta$ is also an endomorphism of $(A, * , [ , , ])$. Then Corollary 3.2 implies that $(A, \tilde{*}, \{ , , \}, \beta)$ is a Hom-LY algebra. \hfill $\square$
\section{Hom-Lie-Yamguti algebras and Hom-Malcev algebras}
In this section we investigate conditions when a Hom-LY algebra reduces to a Hom-Malcev algebra. This consideration is based on the ternary operation $\{,,\}$ of a given Hom-LY algebra $(\mathfrak{m}, [,], \{,,\}, \alpha)$ that could be expressed through its binary one ``$[,]$'' as \\
\par
$ \{ x,y,z \} = -J_{\alpha}(x,y,z) + 2[[x,y], \alpha (z)]$, \hfill (4.1) \\
\\
for all $x,y,z$ in $\mathfrak{m}$, where $J_{\alpha}(x,y,z) := \sigma [[x,y], \alpha (z)]$.
\par
First we recall that a {\it Hom-Malcev algebra} [21] is a Hom-algebra $(A, [,], \alpha )$ such that ``$[,]$'' is skew-symmetric and that the {\it Hom-Malcev identity} \\
\par
$J_{\alpha}(\alpha (x), \alpha (y) ,[x,z]) = [J_{\alpha}(x,y,z),{\alpha}^{2}(x)]$ \hfill (4.2) \\
\\
holds for all $x,y,z$ in $A$. It is observed [21] that when $\alpha = Id$, then (4.2) is the Malcev identity and thus a Hom-Malcev algebra reduces to a Malcev algebra. Other identities, equivalent to the identity (4.2), characterizing Hom-Malcev algebras are found ([21], Proposition 2.8). In [5] it is pointed out another defining identity of Hom-Malcev algebras. This latter identity is the most useful in the proof of the following \\
\par
{\bf Proposition 4.1.} {\it Let $(\mathfrak{m}, [,], \{,,\}, \alpha)$ be a Hom-LY algebra. If its ternary operation ``$\{,,\}$'' expresses through its binary one ``$[,]$'' as in (4.1) for all $x,y,z$ in $\mathfrak{m}$, then $(\mathfrak{m}, [,], \alpha)$ is a Hom-Malcev algebra}. \\
\par
{\bf Proof.} Observe that (4.1) and multiplicativity imply \\
\par
$\{ \alpha (x), \alpha (y), z \} = -J_{\alpha}(\alpha (x), \alpha (y),z) + 2 \alpha ([[x,y],z])$. \hfill (4.3) \\
\\
Then, setting (4.3) in {\bf (B7)}, we get
\par
$-J_{\alpha}(\alpha (x), \alpha (y) ,[u,v]) = [- J_{\alpha}(x,y,u),{\alpha}^{2}(v)] + [{\alpha}^{2}(u), -J_{\alpha}(x,y,v)] $ 
\par
\hspace{4cm} $ + [[2[x,y], \alpha (u)], {\alpha}^{2}(v)] + [{\alpha}^{2}(u), [2[x,y], \alpha (v)]]$
\par
\hspace{4cm}
$ - 2 \alpha ([[x,y],[u,v]])$ \\
and this last equality is written as \\
\par
\par
$J_{\alpha}(\alpha (x), \alpha (y) ,[u,v]) = [J_{\alpha}(x,y,u),{\alpha}^{2}(v)] + [{\alpha}^{2}(u), J_{\alpha}(x,y,v)]$
\par
\hspace{3.5cm} $ -2 J_{\alpha}(\alpha (u), \alpha (v) ,[x,y])$. \hfill (4.4) \\
\\
The expression (4.4) is shown [5] to be equivalent to the Hom-Malcev identity (4.2). Therefore $(\mathfrak{m}, [,], \alpha)$ is a Hom-Malcev algebra. \hfill $\square$ \\
\par
{\bf Remark.} For $\alpha = Id$, the ternary operation (4.1) reduces to the ternary operation, defined by the identity (1.4) in [18], that is considered in Malcev algebras. Thus Proposition 4.1 is the Hom-analogue of the result of K. Yamaguti [18], which is the converse of a result of A.A. Sagle ([13], Proposition 8.3). The Hom-version of the sagle's result is the following \\
\par
{\bf Proposition 4.2.} {\it Let $(\mathfrak{m}, [,], \alpha)$ be a Hom-Malcev algebra and define on $(\mathfrak{m}, [,], \alpha)$ a ternary operation by (4.1). Then} \\
\par
$\{ \alpha (x), \alpha (y),[u,v] \} = [\{ x,y,u \}, {\alpha}^{2}(v)] + [{\alpha}^{2}(u), \{ x,y,v \}] $ \hfill (4.5) \\
\\
{\it for all $u,v,x,y$ in $\mathfrak{m}$}. \\
\par
{\bf Proof.} We write the identity (4.4) as
\par
$-J_{\alpha}(\alpha (x), \alpha (y) ,[u,v]) = [- J_{\alpha}(x,y,u),{\alpha}^{2}(v)] + [{\alpha}^{2}(u), -J_{\alpha}(x,y,v)] $ 
\par
\hspace{4cm}$+ 2 J_{\alpha}(\alpha (u), \alpha (v) ,[x,y])$ \\
i.e.
\par
$-J_{\alpha}(\alpha (x), \alpha (y) ,[u,v]) = [- J_{\alpha}(x,y,u),{\alpha}^{2}(v)] + [{\alpha}^{2}(u), -J_{\alpha}(x,y,v)] $ 
\par
\hspace{4cm}$+ 2 [\alpha ([u,v]),\alpha ([x,y])] + 2[[\alpha (v) ,[x,y]], {\alpha}^{2}(u)]$
\par
\hspace{4cm}$+ 2[[[x,y], \alpha (u)], {\alpha}^{2}(v)]$ \\
or
\par
$-J_{\alpha}(\alpha (x), \alpha (y) ,[u,v]) + 2 [\alpha ([x,y]),\alpha ([u,v])]$
\par
\hspace{4cm}$ = [-J_{\alpha}(x,y,u) + 2[[x,y], \alpha (u)], {\alpha}^{2}(v)] $
\par
\hspace{4cm}$+ [{\alpha}^{2}(u), -J_{\alpha}(x,y,v) + 2[[x,y], \alpha (v)]]$. \\
This last equality (according to (4.1) and using multiplicativity) means that 
\par
$\{ \alpha (x), \alpha (y),[u,v] \} = [\{ x,y,u \}, {\alpha}^{2}(v)] + [{\alpha}^{2}(u), \{ x,y,v \}] $ \\
and therefore the proposition is proved. \hfill $\square$ \\
\par
It could be expected that any Hom-Malcev algebra $(\mathfrak{m}, [,], \alpha)$ with a ternary operation as in (4.1) has a Hom-LY structure. This is for further investigation. Combining Proposition 4.1 and Proposition 4.2, we get the following \\
\par
{\bf Corollary 4.3.} {\it In an anticommutative Hom-algebra $(A, [,], \alpha)$, the Hom-Malcev identity (4.2) is equivalent to (4.5), with $\{ , , \}$ defined by (4.1)}. \hfill $\square$ \\
\par
The untwisted counterpart of Corollary 4.3 is Theorem 1.1 in [18].

\end{document}